\documentclass[11pt]{article}

\usepackage{amsmath,amsfonts,amssymb}
\usepackage{bbm}

\def\|{\vert\vert}

\def\bR {{\mathbb R}}

\def\<{{\langle}}
\def\>{{\rangle}}

\newtheorem{theorem}{Theorem}[section]
\newtheorem{lemma}[theorem]{Lemma}

\newtheorem{proposition}[theorem]{Proposition}

\newtheorem{remark}[theorem]{Remark}

\begin{document}

\title{\vspace*{-1.5cm}
Quadratic Transportation Cost Inequalities Under Uniform Distance For Stochastic Reaction Diffusion Equations Driven by Multiplicative Space-Time White Noise.}

\author{Shijie Shang$^{1}$ and
Tusheng Zhang$^{2,1}$}
\footnotetext[1]{\, School of Mathematics, University of Science and Technology of China, Hefei, China.}
\footnotetext[2]{\, School of Mathematics, University of Manchester, Oxford Road, Manchester
M13 9PL, England, U.K. Email: tusheng.zhang@manchester.ac.uk}
\maketitle

\begin{abstract}
In this paper, we established a  quadratic transportation cost inequality for solutions of stochastic reaction diffusion equations driven by multiplicative space-time white noise based on a new inequality we proved for the moments (under the uniform norm) of the stochastic convolution with respect to space-time white noise, which is of independent interest. The solutions of such stochastic partial differential equations are typically not semimartingales on the state space.
\end{abstract}


\medskip


\section{Introduction}
Let $(X, d)$ be a metric  space with a Borel  probability measure $\mu$. For a measurable subset $A\subset X$ and $r>0$, we denote by $A_r$  the $r$-neighborhood of $A$, namely
$A_r=\{x: d(x,A)<r\}$. We say that $\mu$ has normal concentration on $(X, d)$ if there are constants $C,c>0$ such that for every $r>0$ and every Borel subset $A$ with $\mu(A)\geq \frac{1}{2}$,
\begin{equation}\label{0.1}
1-\mu(A_r)\leq Ce^{-cr^2}.
\end{equation}
 It is well known  that Gaussian measures on $\bR^d$ and uniform measures on the spheres $\mathbb{S}^d$ have normal concentration. In the past decades, many people  established normal concentration properties for various kinds of interesting measures. We  mention the celebrated works of M. Talagrand \cite{T1}, \cite{T2} and \cite{T3}.  We refer the readers to the monograph \cite{L} for a nice exposition of the concentration of measure phenomenon. It turns out that the concentration of measure phenomenon has close connections with entropy and functional inequalities, e.g. Poincare inequalities, logarithmic Sobolev inequalities and transportation cost inequalities. In particular, transportation cost inequalities imply the normal concentration. An elegant, simple proof of this fact is contained in the book \cite{L}. The importance of the topic of the  concentration of measure lies also in its wide applications, e.g. to stochastic finance (see \cite{La}), statistics (see \cite{M}) and the analysis of randomized algorithms (see \cite{DP}).
 \vskip 0.3cm
The concentration of measure for stochastic differential equations and stochastic partial differential equations (SPDEs) has been investigated by many people. Let us  mention several papers which are relevant to our work.  The transportation cost inequalities for stochastic differential equations were obtained by H. Djellout, A. Guillin and L. Wu in \cite{DGW}. The measure concentration for multidimensional diffusion processes with reflecting boundary conditions was considered by S. Pal in \cite{P}. Transportation cost inequalities for solution of stochastic partial differential equations driven by Gaussian noise which is white in time and colored in space were obtained by A. S. Ustunel in \cite{U}. We particularly like to mention the paper \cite{KS} by  D. Khoshnevisan  and A. Sarantsev, which is the starting point of our work. In \cite{KS}, the authors established the quadratic transportation cost inequality under $L^2$-distance for stochastic reaction diffusion equations driven by multiplicative space-time white noise. However, under
the uniform distance they only obtained the quadratic transportation cost inequality for stochastic reaction diffusion equations driven by additive space-time white noise.
As is well known, one of the essential differences between SPDEs driven by colored noise  and SPDEs driven by space-time white noise is that the solution of the later is not a semimartingale and therefore  in particular Ito formula could not be used.
\vskip 0.3cm
The aim of this paper is to prove that under the uniform distance the quadratic transportation cost inequality holds for  stochastic reaction diffusion equations driven by multiplicative space-time white noise. Our new contribution is the $p$th moment inequalities under the uniform norm we obtained for the stochastic convolution with respect to space-time white noise, which is of independent interest.  The significance of the inequality is to allow the order $p$ of the moment to be any positive number, not just for sufficiently large ones. These new estimates allow us to establish the quadratic transportation cost inequality under the uniform norm.
\vskip 0.4cm
The rest of the paper is organized as follows. In Section 2, we recall the notions of measure concentration, transportation cost inequalities and present the framework for stochastic reaction diffusion equations. Section 3 is devoted to the proof of the new moment estimates for stochastic convolutions with respect to space-time white noise under the uniform norm. In Section 4, we prove the quadratic transportation cost inequality.

\section{Preliminaries}\label{S:2}
\setcounter{equation}{0}

In this section, we will recall several results on measure concentration from the monograph \cite{L} and set up the framework of the stochastic reaction diffusion equations driven by space-time white noise.
\vskip 0.4cm
Let $(X,d, \mu)$ be a metric space with a Borel probability measure $\mu$. The concentration function $\alpha_{\mu}(r)$ is defined as
\[\alpha_{\mu}(r):=\sup\left\{1-\mu(A_r): A\subset X, \mu(A)\geq \frac{1}{2}\right\}, \quad r>0.\]
The normal concentration of $\mu$ means that $\alpha_{\mu}(r)\leq Ce^{-cr^2}$ for all $r>0$ with some positive constants $C,c$.

Let $\mu$, $\nu$ be two Borel probability measures on the metric space $(X, d)$. Consider the Wasserstein distance
$$W_2(\nu, \mu):=\left[\inf \int_X\int_X d(x,y)^2\pi(dx,dy)\right]^{\frac{1}{2}}$$
between $\mu$ and $\nu$, where the infimum is taken over all probability measures $\pi$ on the product space $X\times X$ with marginals $\mu$ and $\nu$. Recall that the relative entropy of $\nu$ with respect to $\mu$ is defined by
\[H(\nu|\mu):=\int_X log(\frac{d\nu}{d\mu})d\nu ,\]
if $\nu$ is absolutely continuous with respect to $\mu$, and $+\infty$ if not. We say that the measure $\mu$ satisfies a quadratic transportation cost inequality if there exists a constant $C>0$ such that for all probability measures $\nu$,
\begin{equation}\label{1.2}
 W_2(\nu, \mu)\leq C\sqrt{H(\nu| \mu)}.
 \end{equation}

The following result is taken from \cite{L}.
\begin{proposition}
If $\mu$ satisfies a quadratic transportation cost inequality, then  $\mu$ has normal concentration.
\end{proposition}

\begin{remark} The notion of  concentration of measure phenomenon depends on the underlying  topology of the associated metric space. The stronger the topology, the stronger the concentration.
\end{remark}
\vskip 0.3cm

Before ending this section, let us recall the setup for the stochastic reaction diffusion equations driven by space-time white noise.
Consider the following equation:
\begin{align}\label{3.1}
\left\{
\begin{aligned}
 du(t,x)&=\frac{1}{2}u^{\prime\prime}(t,x)dt+b(u(t,x))dt+ \sigma(u(t,x))W(dt,dx), \quad x\in (0,1), \\
 u(t,0)&=u(t,1)=0  , \quad\quad t> 0,\\
u(0,x)&=u_0(x), \quad x\in (0,1),
\end{aligned}
\right.
\end{align}
where $u_0\in C_0(0,1)$, $W(dt,dx)$ is a space-time white noise on some filtrated probability
space $(\Omega, {\cal F}, {\cal F}_t, P)$, here ${\cal F}_t, t\geq 0$
are the argumented filtration generated by the Brownian
sheet $\{W(t,x); (t,x)\in [0, \infty)\times [0,1]\}$. The coefficients
$b(\cdot), \sigma(\cdot): \mathbb{R}\rightarrow \mathbb{R}$ are deterministic
measurable functions. We say that an adapted,
continuous random field $\{u(t,x): (t,x)\in \mathbb{R}_{+}\times [0,1]\}$ is a solution to the stochastic partial differential equation
(SPDE) (\ref{3.1}) if $t\geq 0$,
\begin{align}\label{3.2}
&\int_0^1u(t,x)\phi(x)dx=\int_0^1u_0(x)\phi(x)dx
+\frac{1}{2}\int_0^tds\int_0^1u(s,x)\phi^{\prime\prime}(x)dx\nonumber\\
&+\int_0^tds\int_0^1b(u(s,x))\phi(x)dx+ \int_0^t\int_0^1\sigma(u(s,x))\phi(x)W(ds,dx),\quad P-a.s.
\end{align}
for any $\phi \in C_0^2(0, 1)$.
It was shown in \cite{W} that $u$ is a solution to SPDE (\ref{3.1}) if and only if $u$ satisfies the following integral equation
\begin{align}\label{3.3}
u(t,x)=&P_tu_0(x)+\int_0^t\int_0^1p_{t-s}(x,y)b(u(s,y))dsdy\nonumber\\
&+ \int_0^t\int_0^1p_{t-s}(x,y)\sigma(u(s,y))W(ds,dy),
\end{align}
where $P_t, t\geq 0$ and $p_{t}(x,y)$ are the corresponding semigroup and the heat kernel associated with the operator $\frac{1}{2}\Delta $ equipped with the Dirichlet boundary condition on the interval $[0, 1]$.
\vskip 0.3cm
Introduce the hypotheses
\begin{itemize}
\item[({\bf H.1})] There exists a constant $L_b$ such that for all $x,y\in \mathbb{R}$,
\begin{align}
|b(x)| \leq\, & L_b(1+|x|),\nonumber\\
|b(x)-b(y)| \leq\, & L_b |x-y|.
\end{align}
 \item[({\bf H.2})] There exist constants $K_{\sigma}$ and $L_{\sigma}$ such that for all $x,y\in \mathbb{R}$,
\begin{align}
|\sigma(x)| \leq\, & K_{\sigma},\nonumber\\
|\sigma(x)-\sigma(y)| \leq\, & L_{\sigma}|x-y|.
\end{align}
\end{itemize}

It is well known (see \cite{W}) that under the hypotheses (H.1) and (H.2), SPDE (\ref{3.1}) admits a unique random field solution $u(t,x)$. In fact, for the existence and uniqueness the diffusion coefficient $\sigma(\cdot)$ needs not to be bounded, the stronger assumption (H.2) is needed for proving the transportation cost inequality.

\section{Moment estimates for stochastic convolution under the uniform norm}
\setcounter{equation}{0}
In this section, we will establish some moment estimates for the stochastic convolution against space-time white noise. Of particular interest are the estimates of the moments of lower order. These bounds will be used later in the paper.
\begin{proposition}\label{estimates 001}
  Let $\{\sigma(s,y): (s,y)\in\mathbb{R}_+\times [0,1]\}$ be a random field such that the stochastic integral against space time white noise is well defined. Then for any $T>0$, $p>10$, there exists a constant $C_{T,p}>0$ such that
  \begin{align}\label{101.1}
    & E\left[\sup_{(t,x)\in[0,T]\times[0,1]}\left|\int_0^t\int_0^1 p_{t-s}(x,y)\sigma(s,y)W(ds,dy) \right|^p\right] \nonumber\\
    \leq & C_{T,p} \int_0^T\sup_{y\in[0,1]}E\left|\sigma(s,y)\right|^p\,ds .
  \end{align}

\end{proposition}

\begin{remark}
  The constant $C_{T,p}$ in (\ref{101.1})
  can be bounded as
\begin{align}
  C_{T,p} < p^{\frac{p}{2}} T^{\frac{p}{4}-\frac{3}{2}} \left(\frac{2}{\pi}\right)^p \left(\frac{1}{\sqrt{2\pi}}\right)^{\frac{p}{2}+1} \left(\frac{6p-8}{p-10}\right)^{\frac{3p}{2}-2} .
\end{align}
\end{remark}

%

\vskip 0.5cm

\noindent {\bf Proof}.
Obviously, we can assume that the right hand side of (\ref{101.1}) is finite.
We employ the factorization method. Choose $\alpha$ such that  $\frac{3}{2p}<\alpha<\frac{1}{4}-\frac{1}{p}$. This is possible because $p>10$. Let
\begin{align}
  (J_{\alpha}\sigma)(s,y):&= \int_0^s\int_0^1 (s-r)^{-\alpha}p_{s-r}(y,z)\sigma(r,z)W(dr,dz), \\
  (J^{\alpha-1}f)(t,x):&= \frac{sin\pi\alpha}{\pi}\int_0^t\int_0^1 (t-s)^{\alpha-1}p_{t-s}(x,y)f(s,y)dsdy.
\end{align}

\noindent By the stochastic Fubini theorem (see Theorem 2.6 in \cite{W}),  for any $(t,x)\in\mathbb{R}_+\times[0,1]$,
\begin{align}\label{103.1}
  \int_0^t\int_0^1 p_{t-s}(x,y)\sigma(s,y)W(ds,dy)=J^{\alpha-1}(J_a\sigma)(t,x).
\end{align}

\noindent Therefore
\begin{align}
  & \sup_{(t,x)\in[0,T]\times[0,1]}\left|\int_0^t\int_0^1 p_{t-s}(x,y)\sigma(s,y)W(ds,dy)\right| \nonumber\\
  =& \sup_{(t,x)\in[0,T]\times[0,1]}\left|J^{\alpha-1}(J_{\alpha}\sigma)(t,x)\right|, \quad P-a.s..
\end{align}
Recall the well-known inequality
\begin{align}
  0\leq p_t(x,y) \leq \frac{1}{\sqrt{2\pi t}}\exp^{-\frac{(x-y)^2}{2t}}, \quad \forall\, x,y \in [0,1].
\end{align}
A straightforward calculation gives
\begin{align}
\label{103.20} \int_0^1 p_t(x,y)\,dy < & 1 , \\
\label{103.2}  \int_0^1 p_t(x,y)^2\,dy =& \sup_{y\in[0,1]}p_t(x,y)\times \int_0^1 p_t(x,y)\,dy \leq C_2 t^{-\frac{1}{2}} , \quad  C_2:=\frac{1}{\sqrt{2\pi}} .
\end{align}
By H\"{o}ler's inquality, (\ref{103.20}) and (\ref{103.2}), we have
{\allowdisplaybreaks\begin{align}\label{104.1}
  & E\sup_{(t,x)\in[0,T]\times[0,1]}\left|\int_0^t\int_0^1 p_{t-s}(x,y)\sigma(s,y)W(ds,dy)\right|^p \nonumber\\
  =& E\sup_{(t,x)\in[0,T]\times[0,1]}\left|\frac{sin\pi\alpha}{\pi} \int_0^t\int_0^1 (t-s)^{\alpha-1} p_{t-s}(x,y)J_{\alpha}\sigma(s,y)\,dsdy\right|^p \nonumber\\
  \leq & \left|\frac{sin\pi\alpha}{\pi}\right|^p E\sup_{(t,x)\in[0,T]\times[0,1]}\bigg\{\int_0^t (t-s)^{\alpha-1} \nonumber\\
  &~~~~~~~~~~~~~\times\left(\int_0^1 p_{t-s}(x,y)|J_{\alpha}\sigma(s,y)|\,dy\right)\,ds\bigg\}^p \nonumber\\
  \leq & \left|\frac{sin\pi\alpha}{\pi}\right|^p E\sup_{(t,x)\in[0,T]\times[0,1]}\Bigg\{\int_0^t (t-s)^{\alpha-1}  \nonumber\\
  &~~~~~~~~~~~~~\times\left(\int_0^1 p_{t-s}(x,y)|J_{\alpha}\sigma(s,y)|^{\frac{p}{2}}\,dy\right)^{\frac{2}{p}}\,ds\Bigg\}^p \nonumber\\
  \leq & \left|\frac{sin\pi\alpha}{\pi}\right|^p E\sup_{(t,x)\in[0,T]\times[0,1]}\Bigg\{\int_0^t (t-s)^{\alpha-1}  \nonumber\\
  &~~~~~~~~~~~~~\times\left(\int_0^1 p_{t-s}(x,y)^2\,dy\right)^{\frac{1}{2}\times\frac{2}{p}}\left(\int_0^1|J_{\alpha}\sigma(s,y)|^p\,dy\right)^{\frac{1}{2}\times\frac{2}{p}}\,ds\Bigg\}^p \nonumber\\
  \leq & \left|\frac{sin\pi\alpha}{\pi}\right|^p C_2 E\sup_{t\in[0,T]}\left\{\int_0^t (t-s)^{\alpha-1-\frac{1}{2p}}\left(\int_0^1|J_{\alpha}\sigma(s,y)|^p\,dy\right)^{\frac{1}{p}}\,ds\right\}^p \nonumber\\
  \leq & \left|\frac{sin\pi\alpha}{\pi}\right|^p C_2 E\sup_{t\in[0,T]}\Bigg[\left(\int_0^t (t-s)^{(\alpha-1-\frac{1}{2p})\frac{p}{p-1}}\,ds\right)^{\frac{p-1}{p}\times p}  \nonumber\\
  &~~~~~~~~~~~~~\times\left(\int_0^t\int_0^1|J_{\alpha}\sigma(s,y)|^p\,dyds\right)^{\frac{1}{p}\times p}\Bigg] \nonumber\\
  \leq & \left|\frac{sin\pi\alpha}{\pi}\right|^p C_2 \times\left(\int_0^T s^{(\alpha-1-\frac{1}{2p})\frac{p}{p-1}}\,ds\right)^{p-1}\times\int_0^T\int_0^1 E|J_{\alpha}\sigma(s,y)|^p\,dyds \nonumber\\
  \leq & C_{T,p}^{\prime}\sup_{(s,y)\in[0,T]\times[0,1]} E\left|\int_0^s\int_0^1 (s-r)^{-\alpha}p_{s-r}(y,z)\sigma(r,z)W(dr,dz)\right|^p ,
\end{align}}
where we have used the condition $\alpha >\frac{3}{2p}$, so that
\begin{align}\label{C_{T,p} prime}
  C_{T,p,\alpha}^{\prime}= & \left|\frac{sin\pi\alpha}{\pi}\right|^p C_2 \times\left(\int_0^T s^{(\alpha-1-\frac{1}{2p})\frac{p}{p-1}}\,ds\right)^{p-1}\times T \nonumber\\
  =& \left|\frac{sin\pi\alpha}{\pi}\right|^p C_2 \left(\frac{p-1}{\alpha p-\frac{3}{2}}\right)^{p-1} T^{\alpha p-\frac{1}{2}} < \infty .
\end{align}

\noindent Applying the BDG inequality (see Proposition 4.4 in \cite{K}) and (\ref{103.2}), we have
\begin{align}
 & \left\Vert\int_0^s\int_0^1 (s-r)^{-\alpha}p_{s-r}(y,z)\sigma(r,z)W(dr,dz)\right\Vert_{L^p(\Omega)}^2 \nonumber\\
 \leq & 4p \int_0^s\int_0^1(s-r)^{-2\alpha}p_{s-r}(y,z)^2\left\Vert\sigma(r,z)\right\Vert_{L^{p}(\Omega)}^2\,drdz \nonumber\\
 \leq & 4p\int_0^s(s-r)^{-2\alpha}\left(\int_0^1 p_{s-r}(y,z)^2\,dz\right)\sup_{z\in[0,1]}\left\Vert\sigma(r,z)\right\Vert_{L^p(\Omega)}^2\,dr  \nonumber\\
  \leq & 4C_2 p \int_0^s(s-r)^{-2\alpha-\frac{1}{2}}\sup_{z\in[0,1]}\left\Vert\sigma(r,z)\right\Vert_{L^p(\Omega)}^2\,dr  \nonumber\\
  \leq & 4C_2 p \left(\int_0^s(s-r)^{(-2\alpha-\frac{1}{2})\times\frac{p}{p-2}}\,dr\right)^{\frac{p-2}{p}} \times\left(\int_0^s \sup_{z\in[0,1]}\left\Vert\sigma(r,z)\right\Vert_{L^p(\Omega)}^p\,dr\right)^{\frac{2}{p}} .
\end{align}
Therefore
\begin{align}\label{104.2}
  & \sup_{(s,y)\in[0,T]\times[0,1]} E\left|\int_0^s\int_0^1 (s-r)^{-\alpha}p_{s-r}(y,z)\sigma(r,z)W(dr,dz)\right|^p \nonumber\\
  \leq & C^{\prime\prime}_{T,p}\times \int_0^T \sup_{z\in[0,1]}E\left|\sigma(r,z)\right|^p\,dr ,
\end{align}
where the condition $\alpha <\frac{1}{4}-\frac{1}{p}$ was used to see that
\begin{align}\label{C_{T,p} prime prime}
  C^{\prime\prime}_{T,p,\alpha}= & (4C_2 p)^{\frac{p}{2}}\times\left(\int_0^T r^{(-2\alpha-\frac{1}{2})\times\frac{p}{p-2}}\,dr\right)^{\frac{p-2}{2}} \nonumber\\
  = & (4C_2 p)^{\frac{p}{2}}\times\left(\frac{p-2}{\frac{p}{2}-2-2\alpha p}\right)^{\frac{p-2}{2}} T^{\frac{p}{4}-1-\alpha p} <\infty .
\end{align}

\noindent Combining (\ref{104.1}) with (\ref{104.2}), we obtain
\begin{align}
  & E\sup_{(t,x)\in[0,T]\times[0,1]}\left|\int_0^t\int_0^1 p_{t-s}(x,y)\sigma(s,y)W(ds,dy)\right|^p \nonumber\\
  \leq & C_{T,p} \int_0^T \sup_{z\in[0,1]}E\left|\sigma(r,z)\right|^p\,dr ,
\end{align}
where
\begin{align}\label{C_{T,p}}
  C_{T,p}= \min_{\frac{3}{2p}<\alpha<\frac{1}{4}-\frac{1}{p}}C^{\prime}_{T,p,\alpha}\times C^{\prime\prime}_{T,p,\alpha} .
\end{align}
In view of (\ref{C_{T,p} prime}), (\ref{C_{T,p} prime prime}) and (\ref{103.2}), a straightforward calculation leads to
\begin{align}
  C_{T,p} < p^{\frac{p}{2}} T^{\frac{p}{4}-\frac{3}{2}} \left(\frac{2}{\pi}\right)^p \left(\frac{1}{\sqrt{2\pi}}\right)^{\frac{p}{2}+1} \left(\frac{6p-8}{p-10}\right)^{\frac{3p}{2}-2} .
\end{align}
%
This completes the proof of the estimate (\ref{101.1}). $\blacksquare$

\vskip 0.5cm
\begin{lemma}\label{estimates 002}
Let $\sigma(s,y)$ be as in Proposition \ref{estimates 001}, then for any $T>0$, $p>10$, $\lambda>0$, there exists a constant $C_{T,p}>0$ such that
  \begin{align}\label{102.1}
    & P\left(\sup_{(t,x)\in[0,T]\times[0,1]}\left|\int_0^t\int_0^1 p_{t-s}(x,y)\sigma(s,y)W(ds,dy) \right|>\lambda\right) \nonumber\\
    \leq & P\left(\int_0^T\sup_{y\in[0,1]}\left|\sigma(s,y)\right|^p\,ds >\lambda^p\right) \nonumber\\
     & +  \frac{C_{T,p}}{\lambda^p}E \min\left\{\lambda^p, \int_0^T\sup_{y\in[0,1]}\left|\sigma(s,y)\right|^p\,ds\right\}.
  \end{align}
  Here the constant $C_{T,p}$ is the same as the constant $C_{T,p}$ in (\ref{101.1}).
\end{lemma}
\vskip 0.3cm
\noindent {\bf Proof}.
For any $\lambda >0$, define
\begin{align}
  \Omega_{\lambda}:=\left\{\omega\in\Omega: \int_0^T \sup_{y\in[0,1]}|\sigma(s,y)|^p\,ds\leq \lambda^p\right\} .
\end{align}
By Chebyshev's inequality, we have
\begin{align}\label{105.1}
    & P\left(\sup_{(t,x)\in[0,T]\times[0,1]}\left|\int_0^t\int_0^1 p_{t-s}(x,y)\sigma(s,y)W(ds,dy) \right|>\lambda\right) \nonumber\\
    \leq & P(\Omega\backslash\Omega_{\lambda}) + P\left(\sup_{(t,x)\in[0,T]\times[0,1]}\left|\int_0^t\int_0^1 p_{t-s}(x,y)\sigma(s,y)W(ds,dy) \right| \mathbbm{1}_{\Omega_{\lambda}} >\lambda\right) \nonumber\\
    \leq & P(\Omega\backslash\Omega_{\lambda}) + \frac{1}{\lambda^p} E\left[\sup_{(t,x)\in[0,T]\times[0,1]} \left|\mathbbm{1}_{\Omega_{\lambda}}\int_0^t\int_0^1 p_{t-s}(x,y)\sigma(s,y)W(ds,dy) \right|^p \right] .
\end{align}
Now, we introduce the random field
\begin{align}
  \widetilde{\sigma}(s,y):= \sigma(s,y)\mathbbm{1}_{\left\{\omega\in\Omega: \ \int_0^s\sup_{y\in[0,1]}|\sigma(r,y)|^p dr \leq\lambda^p\right\}}.
\end{align}
Note that the stochastic integral of $\widetilde{\sigma}(\cdot, \cdot)$ with respect to the space time white noise is well defined.
Since for any $\omega\in\Omega_{\lambda}$,
\begin{align}
  \int_0^t\int_0^1 |\sigma(s,y)-\widetilde{\sigma}(s,y)|^2\,dsdy =0 ,  \quad\forall\,t\in [0,T] ,
\end{align}
by the local property of the stochastic integral (see Lemma \ref{A.1} in Appendix),
\begin{align}
  & \mathbbm{1}_{\Omega_{\lambda}}\int_0^t\int_0^1 p_{t-s}(x,y)\sigma(s,y)W(ds,dy) \nonumber \\
  = & \mathbbm{1}_{\Omega_{\lambda}}\int_0^t\int_0^1 p_{t-s}(x,y)\widetilde{\sigma}(s,y)W(ds,dy), \quad P-a.s..
\end{align}
Hence using the bound (\ref{101.1}), we get
\begin{align}\label{105.2}
  & E\left[\sup_{(t,x)\in[0,T]\times[0,1]} \left|\mathbbm{1}_{\Omega_{\lambda}}\int_0^t\int_0^1 p_{t-s}(x,y)\sigma(s,y)W(ds,dy) \right|^p \right] \nonumber\\
  = & E\left[\sup_{(t,x)\in[0,T]\times[0,1]} \left|\mathbbm{1}_{\Omega_{\lambda}}\int_0^t\int_0^1 p_{t-s}(x,y)\widetilde{\sigma}(s,y)W(ds,dy) \right|^p \right] \nonumber\\
  \leq & E\left[\sup_{(t,x)\in[0,T]\times[0,1]} \left|\int_0^t\int_0^1 p_{t-s}(x,y)\widetilde{\sigma}(s,y)W(ds,dy) \right|^p \right] \nonumber\\
  \leq & C_{T,p} E\int_0^T\sup_{y\in[0,1]}\left|\widetilde{\sigma}(s,y)\right|^p\,ds \nonumber\\
  \leq & C_{T,p} E\min\left\{\lambda^p, \int_0^T\sup_{y\in[0,1]}\left|\sigma(s,y)\right|^p\,ds\right\} .
\end{align}
Combining (\ref{105.1}) with (\ref{105.2}), we obtain (\ref{102.1}).
$\blacksquare$

\begin{proposition}\label{estimates 003}
  Let $\{\sigma(s,y): (s,y)\in\mathbb{R}_+\times [0,1]\}$ be a random field such that the stochastic integral against space time white noise is well defined. Then the following two estimates hold:
\begin{itemize}
  \item [(i)] for any $T>0$, $0<p\leq 10$, $q>10$, there exists a constant $C_{T,p,q}$ such that
  \begin{align}\label{102.2}
    & E\left[\sup_{(t,x)\in[0,T]\times[0,1]}\left|\int_0^t\int_0^1 p_{t-s}(x,y)\sigma(s,y)W(ds,dy) \right|^p\right] \nonumber\\
    \leq & C_{T,p,q} E\left[\int_0^T\sup_{y\in[0,1]}\left|\sigma(s,y)\right|^q\,ds \right]^{\frac{p}{q}} .
  \end{align}

  \item [(ii)] For any $T>0$, $0<p\leq 10$, $\epsilon>0$, there exists a constant $C_{T,p,\epsilon}$ such that
  \begin{align}\label{101.2}
    & E\left[\sup_{(t,x)\in[0,T]\times[0,1]}\left|\int_0^t\int_0^1 p_{t-s}(x,y)\sigma(s,y)W(ds,dy) \right|^p\right] \nonumber\\
    \leq & \epsilon  E\left[\sup_{(s,y)\in[0,T]\times[0,1]}\left|\sigma(s,y) \right|^p\right] + C_{T,p,\epsilon} E\int_0^T\sup_{y\in[0,1]}\left|\sigma(s,y)\right|^p\,ds .
  \end{align}
\end{itemize}
\end{proposition}

\begin{remark}
The significance of the estimates (\ref{102.2}) and (\ref{101.2}) is that they allow $p$ to be small, which is crucial for the proof of the transportation cost inequality in the next section.
\end{remark}

\vskip 0.5cm

\noindent {\bf Proof}.
The estimate (\ref{102.2}) can be easily derived from (\ref{102.1}) and Lemma \ref{A.2} in Appendix as follows:

\begin{align}
  & E\left[\sup_{(t,x)\in[0,T]\times[0,1]}\left|\int_0^t\int_0^1 p_{t-s}(x,y)\sigma(s,y)W(ds,dy) \right|^p\right] \nonumber\\
  = & \int_0^{\infty} p\lambda^{p-1} P\left(\sup_{(t,x)\in[0,T]\times[0,1]}\left|\int_0^t\int_0^1 p_{t-s}(x,y)\sigma(s,y)W(ds,dy) \right|>\lambda\right)d\lambda \nonumber\\
  \leq & \int_0^{\infty} p\lambda^{p-1} P\left(\int_0^T \sup_{y\in[0,1]}|\sigma(s,y)|^q\,ds>\lambda^q\right)d\lambda \nonumber\\
  & + C_{T,p}\int_0^{\infty} p\lambda^{p-1-q}E\min\left\{\lambda^q, \int_0^T \sup_{y\in[0,1]}|\sigma(s,y)|^q\,ds\right\}\,d\lambda \nonumber\\
 = & C_{T,p,q}E \left[\int_0^T\sup_{y\in[0,1]}|\sigma(s,y)|^q\,ds\right]^{\frac{p}{q}} ,
\end{align}
where
\begin{align}
  C_{T,p,q}:=1+ C_{T,p} \frac{q}{q-p},
\end{align}
and the constant $C_{T,p}$ is defined in (\ref{C_{T,p}}).

\vskip 0.3cm

Let us now  prove the assertion (ii) in Proposition \ref{estimates 003}.
From (\ref{102.2}) it follows that for any $q>10$,
{\allowdisplaybreaks\begin{align}\label{106.1}
    & E\left[\sup_{(t,x)\in[0,T]\times[0,1]}\left|\int_0^t\int_0^1 p_{t-s}(x,y)\sigma(s,y)W(ds,dy) \right|^p\right] \nonumber\\
    \leq & C_{T,p,q} E\left[\int_0^T\sup_{y\in[0,1]}\left|\sigma(s,y)\right|^q\,ds \right]^{\frac{p}{q}} \nonumber\\
    \leq & C_{T,p,q} E\left[\sup_{(s,y)\in[0,T]\times[0,1]}|\sigma(s,y)|^{q-p}\times\int_0^T\sup_{y\in[0,1]}\left|\sigma(s,y)\right|^p\,ds \right]^{\frac{p}{q}} \nonumber\\
    = & C_{T,p,q} E\left[\sup_{(s,y)\in[0,T]\times[0,1]}|\sigma(s,y)|^{\frac{(q-p)p}{q}} \times\left(\int_0^T\sup_{y\in[0,1]}\left|\sigma(s,y)\right|^p\,ds\right)^{\frac{p}{q}} \right] \nonumber\\
    \leq & \epsilon E\left[\sup_{(s,y)\in[0,T]\times[0,1]}|\sigma(s,y)|^p\right] + C_{T,p,q}\times C_{T,p,q,\epsilon} E\int_0^T\sup_{y\in[0,1]}\left|\sigma(s,y)\right|^p\,ds  ,
  \end{align}}
where we have used the following Young inequality
\begin{align}
  ab\leq & \frac{\epsilon}{C_{T,p,q}}\,a^{\frac{q}{q-p}} + C_{T,p,q,\epsilon}\,b^{\frac{q}{p}}, \nonumber\\  C_{T,p,q,\epsilon}:= & p\left(\frac{q-p}{\epsilon/C_{T,p,q}}\right)^{\frac{q-p}{p}} q^{-\frac{q}{p}} .
\end{align}
Set
\begin{align}\label{C_{T,p,epsilon}}
  C_{T,p,\epsilon}:=\inf_{q>10} C_{T,p,q}\times C_{T,p,q,\epsilon} .
\end{align}
Now, (\ref{101.2}) follows from (\ref{106.1}) with the constant $C_{T,p,\epsilon}$ defined above.
$\blacksquare$

\section{Quadratic transportation cost inequality}
\setcounter{equation}{0}

In this section, we will show that the law $\mu$ of the random field solution $u(\cdot, \cdot)$ of SPDE (\ref{3.1}), viewed as a probability measure on $C([0, T]\times [0, 1])$, satisfies the quadratic transportation cost inequality, in particular, the normal concentration. First we recall  a lemma proved in \cite{KS} describing the probability measures $\nu$ that are absolutely continuous with respect to $\mu$.
\vskip 0.4cm
Let $\nu\ll \mu$ on $C([0, T]\times [0, 1])$.
Define a new probability measure $Q$ on the filtered probability space $(\Omega, {\cal F}, \{{\cal F}_{t}\}_{0\leq t\leq T}, P)$ by
\begin{align}\label{add 0303.1}
dQ:=\frac{d\nu}{d\mu}(u) dP .
\end{align}
Denote the Radon-Nikodym derivative restricted on ${\cal F}_t$ by
\[M_t:=\left. \frac{dQ}{dP}\right |_{{\cal F}_t}, \quad t\in [0, T].\]
Then $M_t, t\in [0, T]$ forms a $P$-martingale. The following result was proved in \cite{KS}.
\begin{lemma}
There exists an adapted random field $h=\{h(s,x), (s,x)\in [0, T]\times [0,1]\}$ such that $Q-a.s.$ for all $t\in [0, T]$,
\begin{align*}
\int_0^t\int_0^1 h^2(s,x)dsdx<\infty
\end{align*}
and $\widetilde{W}: [0, T]\times [0, 1]\rightarrow \mathbb{R}$ defined by
\begin{align}\label{4.2}
\widetilde{W}(t,x):=W(t,x)-\int_0^t\int_0^x h(s,y)dsdy,
\end{align}
is a Brownian sheet under the measure $Q$. Moreover,
\begin{align}\label{4.3}
M_t=\exp\left (\int_0^t\int_0^1h(s,x)W(ds,dx)-\frac{1}{2}\int_0^t\int_0^1 h^2(s,x)dsdx\right ), \quad Q-a.s.,
\end{align}
and
\begin{align}\label{4.4}
H(\nu|\mu)=\frac{1}{2}E^{Q}\left[\int_0^T\int_0^1 h^2(s,x)dsdx\right],
\end{align}
where $E^{Q}$ stands for the expectation under the measure $Q$.
\end{lemma}

Here is the main result of this section.

\begin{theorem}\label{TCI uniform metric}
Suppose the hypotheses (H.1) and (H.2) hold. Then the law $\mu$ of the solution $u(\cdot, \cdot)$ of SPDE (\ref{3.1}) satisfies the quadratic transportation cost inequality on the space $C([0, T]\times [0, 1])$.  Consequently $\mu$ has normal concentration.
\end{theorem}

\noindent {\bf Proof}.
Take $\nu\ll \mu$ on $C([0, T]\times [0, 1])$.
Define the corresponding measure $Q$ by (\ref{add 0303.1}).
Let $h(t,x)$ be the corresponding random field appeared in Lemma 4.1. Then  the solution $u(t,x)$ of equation (\ref{3.1}) satisfies the following SPDE under the measure $Q$,
 \begin{align}
\label{111.3}  u(t,x)= & P_t u_0(x) + \int_0^t\int_0^1 p_{t-s}(x,y)b(u(s,y))\,dsdy \nonumber\\
  & +\int_0^t\int_0^1 p_{t-s}(x,y)\sigma(u(s,y))\widetilde{W}(ds,dy)  \nonumber\\
  & + \int_0^t\int_0^1 p_{t-s}(x,y)\sigma(u(s,y))h(s,y)\,dsdy .
\end{align}
Consider the solution of the following SPDE:
\begin{align}
\label{111.2}  v(t,x)= & P_t u_0(x) + \int_0^t\int_0^1 p_{t-s}(x,y)b(v(s,y))\,dsdy \nonumber\\
  & +\int_0^t\int_0^1 p_{t-s}(x,y)\sigma(v(s,y))\widetilde{W}(ds,dy).
\end{align}
By Lemma 4.1 it follows that under the measure $Q$, the law of $(v,u)$ forms a coupling of $(\mu, \nu)$. Therefore by the definition of the Wasserstein distance,
\[
W_2(\nu, \mu)^2\leq E^Q\left[\sup_{(t,x)\in[0,T]\times[0,1]}|u(t,x)-v(t,x)|^2\right].
\]
In view of (\ref{4.4}), to prove the quadratic transportation cost inequality
\begin{align}
   W_2(\nu, \mu)\leq \sqrt{2C H(\nu|\mu)} ,
\end{align}
it is sufficient to show that
\begin{align}\label{111.1}
E^{Q}\left[\sup_{(t,x)\in[0,T]\times[0,1]}|v(t,x)-u(t,x)|^2\right]
\leq C E^{Q}\left[\int_0^T\int_0^1 h^2(s,y)\,dsdy\right]
\end{align}
for some independent constant $C$, and assume that the right hand side of (\ref{111.1}) is finite. For simplicity, in the sequel we still denote $E^{Q}$ by the symbol $E$. From (\ref{111.2}) and (\ref{111.3}) it follows that
\begin{align}\label{add 0302.1}
  E \left[\sup_{(t,x)\in[0,T]\times[0,1]}|v(t,x)-u(t,x)|^2\right] \leq 3(I + II +III) ,
\end{align}
where
\begin{align*}
  I := & E\left[\sup_{(t,x)\in[0,T]\times[0,1]}\left|\int_0^t\int_0^1 p_{t-s}(x,y)\big[b(v(s,y))-b(u(s,y))\big]\,dsdy\right|^2\right] , \nonumber\\
  II := & E\left[\sup_{(t,x)\in[0,T]\times[0,1]}\left|\int_0^t\int_0^1 p_{t-s}(x,y)\big[\sigma(v(s,y))-\sigma(u(s,y))\big]\tilde{W}(ds,dy)\right|^2\right] , \nonumber\\
  III := & E\left[\sup_{(t,x)\in[0,T]\times[0,1]}\left|\int_0^t\int_0^1 p_{t-s}(x,y)\sigma(u(s,y))h(s,y)\,dsdy\right|^2\right] .
\end{align*}

\noindent By Holder's inequality and (\ref{103.2}), the term $I$ can be estimated as follows:
{\allowdisplaybreaks\begin{align}\label{term I}
  I \leq & L_b^2 E\left[\sup_{(t,x)\in[0,T]\times[0,1]}\left|\int_0^t\int_0^1 p_{t-s}(x,y)|v(s,y)-u(s,y)|\,dsdy\right|^2\right]  \nonumber\\
  \leq & L_b^2 E\Bigg\{\sup_{(t,x)\in[0,T]\times[0,1]}\bigg[\left(\int_0^t\int_0^1 p_{t-s}(x,y)^2\,dsdy\right) \nonumber\\
  &~~~~~~~~~~~~~~~~~~~~~~~~~\times\left(\int_0^t\int_0^1 |v(s,y)-u(s,y)|^2\,dsdy\right)\bigg]\Bigg\}  \nonumber\\
  \leq & \sqrt{\frac{2T}{\pi}}L_b^2 E\int_0^T\int_0^1 |v(s,y)-u(s,y)|^2\,dsdy  \nonumber\\
  \leq & \sqrt{\frac{2T}{\pi}}L_b^2 \int_0^T  E\left[\sup_{(r,y)\in[0,s] \times[0,1]}|v(r,y)-u(r,y)|^2\right]\,ds .
\end{align}}

\noindent For the term $II$, applying the estimate (\ref{101.2}) we obtain that for any $\epsilon>0$,
\begin{align}\label{term II}
  II \leq & \epsilon E\left[\sup_{(t,x)\in[0,T]\times[0,1]}|\sigma(v(t,x))-\sigma(u(t,x))|^2\right] \nonumber\\
  & + C_{T,2,\epsilon} E\int_0^T\sup_{y\in[0,1]}\left|\sigma(v(s,y))-\sigma(u(s,y))\right|^2\,ds \nonumber\\
  \leq & \epsilon L_{\sigma}^2 E\left[\sup_{(t,x)\in[0,T]\times[0,1]}|v(t,x)-u(t,x)|^2\right] \nonumber\\
  & + C_{T,2,\epsilon}L_{\sigma}^2 \int_0^T E\left[\sup_{(r,y)\in[0,s]\times[0,1]}\left|v(r,y)-u(r,y)\right|^2\right]\,ds .
\end{align}

\noindent The term $III$ can be bounded as follows:
\begin{align}\label{term III}
  III \leq & K_{\sigma}^2 E\Bigg\{\sup_{(t,x)\in[0,T]\times[0,1]}\bigg[\left(\int_0^t\int_0^1 p_{t-s}(x,y)^2\,dsdy\right) \nonumber\\
  &~~~~~~~~~~~~~~~~~~~~~~~~~~\times\left(\int_0^t\int_0^1 h^2(s,y)\,dsdy\right)\bigg]\Bigg\} \nonumber\\
  \leq & \sqrt{\frac{2T}{\pi}} K_{\sigma}^2 E\left[\int_0^T\int_0^1 h^2(s,y)\,dsdy\right].
\end{align}
Set
\begin{align}
  Y(t):= E\left[\sup_{(s,x)\in[0,t]\times[0,1]}|v(s,x)-u(s,x)|^2\right] .
\end{align}
Putting (\ref{add 0302.1})-(\ref{term III}) together, we obtain
\begin{align}\label{4.1}
  Y(T) \leq & 3\sqrt{\frac{2T}{\pi}}L_b^2 \int_0^T Y(s)\,ds + 3 \epsilon L_{\sigma}^2 Y(T) + 3 C_{T,2,\epsilon}L_{\sigma}^2 \int_0^T Y(s)\,ds \nonumber\\
  & + 3\sqrt{\frac{2T}{\pi}} K_{\sigma}^2 E\left [\int_0^T\int_0^1 h^2(s,y)\,dsdy\right ].
\end{align}
Recall that(see e.g. Theorem 3.13 in \cite{DKZ})
\begin{align}
  & E\left[\sup_{(t,x)\in[0,T]\times[0,1]}|u(t,x)|^2\right]< \infty , \\
  & E\left[\sup_{(t,x)\in[0,T]\times[0,1]}|v(t,x)|^2\right]< \infty .
\end{align}
Hence $Y(T)<\infty$ for any $T>0$.
Taking any $\epsilon<\frac{1}{3 L_{\sigma}^2}$, we deduce from (\ref{4.1}) that
\begin{align}\label{add 0308.1}
    Y(T) \leq & \frac{3L_b^2}{1-3 \epsilon L_{\sigma}^2}\sqrt{\frac{2T}{\pi}} \int_0^T Y(s)\,ds + \frac{3 C_{T,2,\epsilon} L_{\sigma}^2}{1-3\epsilon L_{\sigma}^2} \int_0^T Y(s)\,ds \nonumber\\
  & + \frac{3K_{\sigma}^2}{1-3\epsilon L_{\sigma}^2}\sqrt{\frac{2T}{\pi}} E\left[\int_0^T\int_0^1 h^2(s,y)\,dsdy\right] .
\end{align}
Clearly, (\ref{add 0308.1}) still holds if we replace $T$ with any $t\in[0,T]$. Applying Gronwall's inequality, we obtain
\begin{align}
  Y(T) \leq & K_{\sigma}^2\inf_{0<\epsilon<\frac{1}{3L_{\sigma}^{2}}}\left\{\frac{3}{1-3\epsilon L_{\sigma}^2}\sqrt{\frac{2T}{\pi}} \exp\left(\frac{3L_b^2T}{1-3\epsilon L_{\sigma}^2}\sqrt{\frac{2T}{\pi}} +\frac{3 C_{T,2,\epsilon} L_{\sigma}^2T}{1-3\epsilon L_{\sigma}^2}\right)\right\} \nonumber\\
  & \times E\left[\int_0^T\int_0^1 h^2(s,y)\,dsdy\right],
\end{align}
where the constant $C_{T,2,\epsilon}$ is defined in (\ref{C_{T,p,epsilon}}) with $p=2$. This proves (\ref{111.1}), hence completes the proof of Theorem \ref{TCI uniform metric}. $\blacksquare$

\section{Appendix}
\setcounter{equation}{0}

The following local property of the Walsh stochastic integral against space-time white noise is similar to that of the Ito integral.
\begin{lemma}\label{A.1}
Let $\{\sigma(t,x): (t,x)\in [0,T]\times [0,1]\}$ be a random field such that the stochastic  integral against space time white noise is well defined. Let $\Omega_0 \subset \Omega$ be a measurable subset such that for a.s. $\omega\in\Omega_0$,
\begin{align}
  \int_0^T\int_0^1 |\sigma(t,x)|^2\,dtdy=0.
\end{align}
Then for a.s. $\omega\in\Omega_0$,
\begin{align}
  \int_0^T\int_0^1 \sigma(t,x) W(dt,dx) =0 .
\end{align}

\end{lemma}

\noindent {\bf Proof}. The local property  can be similarly proved as that of Ito integral. We only outline the proof here. Firstly, we note that the local property obviously holds when $\sigma(\cdot, \cdot)$ is a simple process. When $\sigma(t, x)$ is a bounded, continuous random field, we can prove the local property through an approximation of  $\sigma$ by a sequence of simple processes. For the general random field $\sigma(\cdot, \cdot)$, the local property can be proved by further two approximations, first by  bounded random fields and then by  continuous random fields.
$\blacksquare$

\begin{lemma}\label{A.2}
  Let $X\geq 0$ be a random variable, then for any $0<p<q$,
  \begin{gather}
\label{Appendix.1}    EX^p=\int_0^{\infty} px^{p-1}P(X>x)\,dx , \\
\label{Appendix.2}    \int_0^{\infty} \frac{E\min\{x^q, X\}}{x^q} px^{p-1}\,dx = \frac{q}{q-p} E\left[X^{\frac{p}{q}}\right] .
  \end{gather}

\end{lemma}

\noindent {\bf Proof}. (\ref{Appendix.1}) and (\ref{Appendix.2}) can be easily proved by Fubini theorem. (\ref{Appendix.2}) is similar to Lemma 2 in \cite{I}, for completeness, we provide the proof here.
\begin{align}
      \int_0^{\infty} \frac{E\min\{x^q, X\}}{x^q} px^{p-1}\,dx = & E\int_0^{X^{\frac{1}{q}}} px^{p-1}\,dx + E\left[X\int_{X^{\frac{1}{q}}}^{\infty} p x^{p-1-q}\,dx\right] \nonumber\\
      = & E\left[X^{\frac{p}{q}}\right] - \frac{p}{p-q} E\left[X\left(X^{\frac{1}{q}}\right)^{p-q}\right] \nonumber\\
      = & \frac{q}{q-p} E\left[X^{\frac{p}{q}}\right] .
\end{align}
$\blacksquare$

\vskip 0.3cm
\noindent{\bf Acknowledgement}.  This work is partially supported by NNSF of China (11671372,  11431014, 11721101).

\end{document}